\documentclass[11pt,fleqn,twoside]{article}
\usepackage{amsfonts,amssymb,latexsym}
\makeatletter
\newcommand{\prava}[1]{\small\it
\begin{flushleft}
Copyright \copyright \ 1999 by  #1
\end{flushleft}}

\newcommand{\name}[1]{\begin{flushleft}
                       \LARGE \bf #1
                       \end{flushleft}\vspace{-3mm}}

\newcommand{\Author}[1]{\begin{flushleft}
                       \it #1 \end{flushleft}}

\newcommand{\Adress}[1]{\begin{flushleft}
                       \it #1 \end{flushleft}}

\newcommand{\Date}[1]{\begin{flushleft}
                      \small  \it #1 \end{flushleft}}

\newcommand{\ehkol}{Author \ name}
\newcommand{\ohkol}{Article \ name}
\renewcommand{\@evenhead}{
\hspace*{-3pt}\raisebox{-15pt}[\headheight][0pt]{\vbox{\hbox to \textwidth 
{\thepage \hfil \ehkol}\vskip4pt \hrule}}}
\renewcommand{\@oddhead}{
\hspace*{-3pt}\raisebox{-15pt}[\headheight][0pt]{\vbox{\hbox to \textwidth 
{\ohkol \hfil \thepage}\vskip4pt\hrule}}}
\renewcommand{\@evenfoot}{}
\renewcommand{\@oddfoot}{}

\newcommand{\be}{\begin{equation}}
\newcommand{\ee}{\end{equation}}
\newcommand{\ba}{\hspace*{-5pt}\begin{array}}
\newcommand{\ea}{\end{array}}

\newcommand{\ds}{\displaystyle}
\makeatother

\begin{document}

\thispagestyle{empty}
\setcounter{page}{269}
\renewcommand{\ehkol}{B.A. Kupershmidt}
\renewcommand{\ohkol}{Remarks on Quantization of Classical $r$-Matrices}

\begin{flushleft}
\footnotesize \sf
Journal of Nonlinear Mathematical Physics \qquad 1999, V.6, N~3,
\pageref{kupershmidt_6-fp}--\pageref{kupershmidt_6-lp}.
\hfill {\sc Letter}
\end{flushleft}

\vspace{-5mm}

\renewcommand{\footnoterule}{}
{\renewcommand{\thefootnote}{}
 \footnote{\prava{B.A. Kupershmidt}}}

\name{Remarks on Quantization of Classical 
{\mathversion{bold}$r$}-Matrices}\label{kupershmidt_6-fp}

\Author{Boris A. KUPERSHMIDT}

\Adress{Department of Mathematics, University of Tennessee
Space Institute,\\ 
Tullahoma, TN 37388,  USA\\
E-mail: bkupersh@utsi.edu}

\Date{Received May 6, 1999; Revised May 26, 1999; Accepted
June 1, 1999}

\begin{abstract}
\noindent
If a classical $r$-matrix $r$ is skewsymmetric, its quantization $R$ can 
lose the skewsymmetry property.  Even when $R$ is skewsymmetric, it may 
not be unique.
\end{abstract}

\noindent
Let $r$ be a classical $r$-matrix.  In general, it means that we have 
a family of vector spaces $\{V_\alpha\}$,  $\alpha   \in  {\mathcal{A}}$, and a 
collection of linear operators
\be
r (\alpha, \beta): \ V_\alpha \otimes V_\beta \rightarrow V_\beta \otimes 
V_\alpha, \qquad  \forall \; \alpha \not= \beta \in {\mathcal{A}}, 
\ee
satisfying the misnamed ``Classical Yang-Baxter'' equation (CYB)
\be
\ba{l}
\ds [c(r)] ^{\varphi \psi \xi}_{ijk} (\alpha, \beta, \gamma):
=\left (r (\alpha, 
\beta)_{ij}^{s \varphi} r(\beta, \gamma)_{sk}^{\xi \psi} + r (\alpha, \beta)
_{ij}^{\psi s} r(\alpha, \gamma)^{\xi \varphi}_{sk}\right)  
\vspace{2mm}\\
\phantom{[c(r)] ^{\varphi \psi \xi}_{ijk} (\alpha, \beta, \gamma): =}
+ c.p. (i, j, k; \varphi, \psi, \xi; \alpha, \beta, \gamma) = 0 
\ea
\ee
where ``$c.p$'' stands for the sum on cyclically permuted triples of indices 
indicated, and $r (\alpha, \beta)^{uv}_{ij}$ are the matrix elements of the 
operators $r(\alpha, \beta)$ (1) in a collection of fixed basises:
\be
r(\alpha, \beta) \left(e^\alpha_{i} \otimes e^\beta_j\right) = r(\alpha, \beta)^{\ell 
k}_{ij} e^\beta_{\ell} \otimes e^{\alpha}_k; 
\ee
the convention of summation over repeated upper-lower indices is in force.

In most applications, all the vector spaces $V_\alpha$ are isomorphic to each 
other, $V_\alpha \approx V$; in addition, often, -- but not always, -- the 
operator $r: V \otimes V \rightarrow V \otimes V$ is skewsymmetric:
\be
P r P = - r, \qquad r^{k \ell}_{ij} = - r^{\ell k}_{ji}, 
\ee
where $P$ is the permutation operator,
\be
P (x \otimes y) = y \otimes x. 
\ee
We shall consider this particular framework from now on.

To quantize a given $r$-matrix $r$ is to find an operator family
\be
R = R (h): \ V \otimes V \rightarrow V \otimes V, 
\ee
depending upon a parameter $h$, such that
\be
R (h) = P + h r + O (h^2), 
\ee
and $R$ satisfies the Artin braid relation (also misnamed as the 
``Quantum Yang-Baxter'' equation, QYB):
\be
R^{12} R^{23} R^{12} = R^{23} R^{12} R^{23}, 
\ee
where this equality of operators acting on $V \otimes V \otimes V$ employs
the standard notation
\be
R^{12} (x \otimes y \otimes z) = R(x \otimes y) \otimes z, \qquad
 R^{23}  (x \otimes y \otimes z) = x \otimes R(y \otimes z). 
\ee

How does the skewsymmetry condition on $r$, (4), translate into $R = R(h)$?

There are at least two possible, {\it logically  independent}, answers, only one of 
which is correct.

The first one is what is commonly accepted in the literature under the name 
of ``unitarity'':
\renewcommand{\theequation}{\arabic{equation}{\rm a}}
\setcounter{equation}{9}
\be
R(h)^{-1} = R (h), 
\ee
or 
\renewcommand{\theequation}{\arabic{equation}{\rm b}}
\setcounter{equation}{9}
\be
R(q)^{-1} = R(q), 
\ee
in the multiplicative notation $q = e^h$.

The second one I shall call, for want of a better term, the mirror symmetry:
\renewcommand{\theequation}{\arabic{equation}}
\setcounter{equation}{10}
\be
R^{{\mathcal{M}}} (h) = R(-h).
\ee
Here $R^{{\mathcal{M}}}$ is the operator acting as the mirror image of $R$.  If
\be
R\left(e_i \otimes e^\prime_j\right) = R^{k \ell}_{ij} e^\prime_k \otimes e_\ell, 
\ee
then
\be
R^{{\cal{M}}} \left(R^{k \ell}_{ij} e_\ell \otimes e^\prime_\kappa\right) = e^\prime_j \otimes e_i. 
\ee
This definition, useful as it is, is {\it not} connected to skewsymmetry of $r$.

The classical $r$-matrix $r$ appears as the $h^1$-term in the $h$-expansion 
of the Quantum $R$-matrix $R(h)$ around $h=0$.  The terms in $h$ of orders 
higher than~1 recede away in the quasiclassical passage.  The examples that 
follow demonstrate that these higher-order terms can have distinctly 
anti-Prussian character and break out strict orders and symmetries.  
(In~[1] Drinfel'd proved that every skewsymmetric classical $r$-matrix $r$ 
represents $h^1$-part of some skewsymmetric Quantum $R-$matrix $R$.  The 
question of additional parameters in $R$ was not addressed there, or 
elsewhere.)

In the $1^{st}$ example, $\mbox{dim}\, (V) =2$ and the $R$-matrix $R=R(h;\theta)$ 
acts on $V \otimes V$ (in a chosen basis) as
\be
R \left(e_0 \otimes e^\prime_0\right) = e^\prime_0 \otimes e_0, 
\ee
\be
R \left(e_0 \otimes e^\prime_1\right) = \left(e^\prime_1 + h e^\prime_0\right) \otimes e_0, 
\ee
\be
R \left(e_1 \otimes e^\prime_0\right) = e^\prime_0 \otimes (e_1 - h e_0), 
\ee
\be
R \left(e_1 \otimes e^\prime_1\right) = e^\prime_1 \otimes e_1 + \theta
 h^2 e^\prime_0 \otimes e_0. 
\ee
Here $\theta$ is an arbitrary constant.  The Artin relation (8) is easily 
verified.  The $h^1$-terms comprise the $r$-matrix
\be
r^{k\ell}_{ij} = \delta^k_0 \delta^\ell_0 \left(\delta^{01}_{ij} - \delta^{10}_{ij} \right) 
\ee
which is obviously skewsymmetric.  The $R$-matrix $R(h; \theta)$ is, however, 
not unitary unless $\theta = 0$.  Also, it's easy to see that
\be
R^{{\mathcal{M}}} (h; \theta) = R(-h; - \theta), 
\ee
so that this $R$-matrix is not mirror-symmetric either, again unless
$\theta = 0$.

Our 2$^{nd}$ example is a little bit more elaborate, with $\mbox{dim}\,(V) =3$.  
Here the $R$-matrix {\it is} both skewsymmetric and mirror-symmetric, but it 
depends upon one extra parameter, in addition to the quantization parameter 
$h$, thus exhibiting clearly nonuniqueness of quantization of classical 
$r$-matrices.

Fixing a basis $(e_0, e_1, e_2)$ in $V$, we set
\renewcommand{\theequation}{\arabic{equation}.1}
\setcounter{equation}{19}
\be
R \left(e_0 \otimes e^\prime_0\right) = e^\prime_0 \otimes e_0, 
\ee
\renewcommand{\theequation}{\arabic{equation}.2}
\setcounter{equation}{19}
\be
R \left(e_0 \otimes e^\prime_1\right) = \left(e^\prime_1 +h e^\prime_0\right) \otimes e_0, 
\ee
\renewcommand{\theequation}{\arabic{equation}.3}
\setcounter{equation}{19}
\be
R \left(e_1 \otimes e^\prime_0\right) = e^\prime_0 \otimes (e_1 - h e_0), 
\ee
\renewcommand{\theequation}{\arabic{equation}.4}
\setcounter{equation}{19}
\be
R \left(e_1 \otimes e^\prime_1\right) = e^\prime_1 \otimes e_1;
\ee
\renewcommand{\theequation}{\arabic{equation}.1}
\setcounter{equation}{20}
\be
R \left(e_0 \otimes e^\prime_2\right) = \left(e^\prime_2 + h e^\prime_1 + {h^2 \over 
2} e^\prime_0\right)  \otimes e_0,  
\ee
\renewcommand{\theequation}{\arabic{equation}.2}
\setcounter{equation}{20}
\be
R \left(e_2 \otimes e^\prime_0\right) = e^\prime_0 \otimes \left(e_2 - h e_1 + {h ^2 
\over 2} e_0\right), 
\ee
\renewcommand{\theequation}{\arabic{equation}.3}
\setcounter{equation}{20}
\be
R \left(e_1 \otimes e^\prime_2\right) = e^\prime_2 \otimes (e_1 + h e_0) + h^2 
\left({1 \over 2} e^\prime_1 + \lambda h e^\prime_0\right) \otimes e_0, 
\ee
\renewcommand{\theequation}{\arabic{equation}.4}
\setcounter{equation}{20}
\be
R \left(e_2 \otimes e^\prime_1\right) = \left(e^\prime_1 -h e^\prime_0 \right)\otimes e_2 
+h^2 e^\prime_0 \otimes \left({1 \over 2} e_1 - \lambda h e_0\right), 
\ee
\renewcommand{\theequation}{\arabic{equation}.5}
\setcounter{equation}{20}
\be
\ba{l}
\ds R \left(e_2 \otimes e^\prime_2\right) = e^\prime_2 \otimes \left(e_2 + h e_1 + {h ^2 
\over 2} e_0\right) 
\vspace{3mm}\\
\ds \phantom{R \left(e_2 \otimes e^\prime_2\right) = }
- h e^\prime_1 \otimes \left(e_2 + {\tilde \lambda} h^2 
e_0\right) + h^2 e^\prime_0 \otimes \left({1 \over 2} e_2 + {\tilde \lambda} h e_1\right). 
\ea
\ee
Here $\lambda$ is the new free parameter, and
\renewcommand{\theequation}{\arabic{equation}}
\setcounter{equation}{21}
\be
{\tilde \lambda} = \lambda - {1 \over 4}. 
\ee

From formulae (20) we see that the previous example (14)--(17) is embedded 
into this one, with $\theta =0$.  It's immediate to check that 
\be
R (h; \lambda)^2 = {\mathbf{1}}, 
\ee
\be
R^{{\mathcal{M}}} (h; \lambda) = R (-h; \lambda), 
\ee
so that our $R$-matrix is both skewsymmetric and mirror-symmetric.  Also, the 
$h^1$-part of $R(h; \lambda)$ is given by the flag-type formula
\be
r^{k\ell}_{ij} = (i - c)\delta^\ell_i \delta^k_{j-1} - (j -c) 
\delta^\ell_{i-1}  \delta^k_j, \qquad 0 \leq i, j, k, \ell \leq  \mbox{dim}\, (V)-1, 
\ee
where $c$ is an arbitrary constant.  [In our case $c=1$, but this constant 
can be adjusted to any desired value by an appropriate nonlinear transformation; 
in particular, we can make 
\be
c = {\mbox{dim}\, (V) - 1 \over 2} 
\ee
to have the determinant in $GL(V)$ being central in the induced Lie-Poisson 
structure~[2].  In this language, the $R$-matrix (20)--(21) defines the Quantum 
Group $\mbox{Mat}_{h;\lambda}(3)$, a 3-di\-mensional analog of the 2-dimensional 
Quantum Group $\mbox{Mat}_{h}$(2).]  The checking of the Artin relation for the 
$R$-matrix (20)--(21) is easy but tedious; the mirror property (24) cuts 
the verification procedure by $1/3$; there are still more symmetries present in 
this $R$-matrix which will allow another $1/3$ of the checking labor to 
be avoided.

How many additional constants should one expect when quantizing a 
skewsymmetric classical $r$-matrix and requiring the Quantum $R$-matrix 
to be skewsymmetric and mirror-symmetric?  For the case of the $r$-matrix (25), I expect 
the total number of additional parameters (the $\lambda$'s) to be
\be
\mbox{dim}\, (V) -2, 
\ee
and in general it probably could never be larger no matter what $r$ is; 
dropping off the mirror-symmetry condition increases the number of possible 
parameters by~1.

\label{kupershmidt_6-lp}

\end{document}